\documentclass{article}
\usepackage{amssymb}

\usepackage{graphicx}
\usepackage{amsmath}

\newtheorem{theorem}{Theorem}

\newtheorem{definition}[theorem]{Definition}

\newtheorem{lemma}[theorem]{Lemma}

\begin{document}

\title{The It\^o transform for a general class of pseudo-differential operators}
\author{R\'emi L\'eandre   \\
Institut de Math\'ematiques. Universit\'e de Bourgogne.
\\ 21000. Dijon. FRANCE.
} \maketitle

\begin{abstract}
We give an It\^o formula for a general class of pseudo-differential operators.
\end{abstract} 
\section{Introduction}
Let us recall what is the It\^o formula for a purely discontinuous martingale $t \rightarrow M_t$ with values in $\mathbb{R}$ [1]. Let $f$ be a $C^2$ function on $\mathbb{R}$. We have
\begin{equation}f(M_t) = f(M_0) + \int_0^tf'(M_{s-})\delta M_s + \sum_{s\leq t}f(M_s)-f(M_{s-})-f'(M_{s-})\Delta M_s\end{equation}
It is the generalization of the celebrated It\^o formula for the Brownian motion $t \rightarrow B_t$ on $\mathbb{R}$ [1]
\begin{equation}f(B_t) = f(B_0) + \int_0^tf'(B_s)\delta B_s + 1/2\int_0^tf"(B_s)ds\end{equation}
A lot of of stochastic analysis tools for diffusions were translated by L\'eandre in semi-group theory in [6], [7], [8], [10], [11], 
 [14], [15], [16], [18]. Some basical tools of stochastic analysis for the study of jump processes
 were translated by L\'eandre in semi-group theory in [11], [12], [19]. For review on that, we refer to the review of L\'eandre [9], [17].

L\'eandre has extended the It\^o formula for the Brownian motion to the case of some classical partial differential equations in [19], [21], [22], [23]. 
In such a case, there is until now no convenient
measure on a convenient path space associated to this partial differential equation.
In [23], we have extended the It\^o formula for jump process for an integro-differential generator when
 there is until now no stochastic process associated. Jump processes are generically generated by pseudo-differential operators which satisfy the maximum 
principle [5].

In this paper, we give an It\^o formula for a general class of positive elliptic pseudo differential operators. 
For material on pseudo-differential operators, we refer on [2], [3], [4] and [5]. Since the considerations below on pseudo-differential operators
 are more and less classical, we won't enter in the technical details of the 
proof.

\section{The two semi-groups}
Let $\hat{u}$ be the Fourier transform of a smooth function $u$ on $\mathbb{R}^d$. Let $a(x, \xi)$ be a global symbol of 
order $m$ on $\mathbb{R}^d$. It is a smooth function from $\mathbb{R}^d \times \mathbb{R}^d$ into $\mathbb{C}$ such that for all $k$, $k'$
\begin{equation}\sup_{x \in \mathbb{R}^d}\vert D^k_xD^{k'}_\xi a(x, \xi)\vert \leq C_{k,k'}\vert \xi \vert^{m-k'}\end{equation}
We say that a global symbol of order $m$ is elliptic if for $\vert \xi \vert >M$
\begin{equation}\inf_{x \in \mathbb{R}^d}\vert a(x, \xi)\vert \geq C_M\vert \xi \vert^m \end{equation}
We consider the proper pseudodifferential operator associated to the symbol $a$: the Fourier transform of $L_0u$ is given by
\begin{equation} \int_{\mathbb{R}^d}a(x, \xi)\hat{u}(\xi)d\xi \end{equation}
We consider its adjoint $L_0^*$ on $L^2(dx)$ and we put $L = L_0^*L_0$.

All the considerations of [2] which were valid on a compact subset of $\mathbb{R}^d$ are still true because (3) and (4) are valid globally.
In particular, $L$ is essentially selfadjoint on $L^2(dx)$ and generates a contraction semi-group $P_t$ on $L^2(dx)$.

Let us consider a smooth function $f$ from $\mathbb{R}^d$ into $\mathbb{R}$
 with compact support and a smooth function $v$ with compact support from $\mathbb{R}^d \times \mathbb{R}$
 into $\mathbb{C}$. $(x,y)$ denotes the generic element of $\mathbb{R}^d \times \mathbb{R}$. We consider the smooth function 
 from $\mathbb{R}^d$ into $\mathbb{R}$ $\hat{v}$
\begin{equation} \hat{v}(x) = v(x, f(x))\end{equation}
We consider the function $\overline{v}$ from $\mathbb{R}^d \times \mathbb{R}$ into $\mathbb{C}$
\begin{equation}(x,y) \rightarrow v(x, y + f(x))\end{equation}
We apply $L$to $\overline{v}$, $y$ being frozen. We get a function $L\overline{v}$. We put
\begin{equation}(\hat{L}v)(x,y) = (L\overline{v})(x, y-f(x))\end{equation}
\begin{definition}$\hat{L}$ is called the It\^o transform of $L$.\end{definition}
We remark that $(x,y) \rightarrow (x, y+f(x))$ is a diffeomorphism of
$\mathbb{R}^d \times \mathbb{R}$ which keeps the measure $dx \otimes dy$ invariant. This shows:
\begin{theorem}$\hat{L}$ is positive symmetric on $L^2(dx \otimes dy)$. It admits therefore a self-adjoint extension still denoted $\hat{L}$. This self-adjoint 
extension generates a semigroup $\hat{P}_t$ of contraction on $L^2(dx \otimes dy)$ \end{theorem}
We get
\begin{theorem}(It\^o formula)We have the relation for all smooth function $v$ with compact support
\begin{equation}P_t(\hat{v})(x) = (\hat{P}_t(v))(x, f(x))\end{equation}\end{theorem}
{\bf{Remark:}}If we consider the generator $L = \sum X_i^2$ where the $X_i$ are smooth vector fields, $\hat{L} = \sum \hat{X}_i^2$ where 
\begin{equation}\hat{X}_i = (X_i, <X_i, df>)\end{equation}
which corresponds to the generator of [19], [21], [22]. Analogous remark holds for the considerations of [23].

\section{Proof of the It\^o formula}
\begin{lemma}If $v$ is a smooth function on $\mathbb{R}^d \times \mathbb{R}$ whose all derivatives belong to $L^2$, $\hat{P}_tv$ is still a smooth function whose all
 derivatives belong to $L^2$.\end{lemma}
{\bf{Proof:}} Let 
\begin{equation}\overline{L}_1 = \hat{L} + (-{\partial^2 \over \partial y^2})^{m/2}\end{equation}
$\overline{L}$ commute with $\hat{L}$. Therefore, for all $k$
\begin{equation}(\overline{L}_1^k)\hat{P}_t = (\hat{P}_t)(\overline{L}_1^k)\end{equation}
If $v$ satisfies the hypothesis, $\hat{P}_tv$ belongs to the domain of $\overline{L}_1^k$. But $\overline{L}$ is the transform of 
\begin{equation}\tilde{L} = {L} + (-{\partial^2 \over \partial y^2})^{m/2}\end{equation}
under the change of variable $(x,y) \rightarrow (x, y+f(x))$. Therefore $\hat{P}_tv$ belongs to the domain of $\tilde{L}^k$. The result arises by 
Garding inequality.$\diamondsuit$

Let $\phi$ be a smooth function from $\mathbb{R}^d$ into $[0,1]$, equals to 0 if $\vert \xi \vert \geq 2$ and equals to 1 if $\vert \xi \vert \leq 1$. We consider the global symbol
\begin{equation}a_\lambda(x,\xi) = \phi(\xi/\lambda)a(x, \xi)\end{equation}
and the operator $L_{0, \lambda}$, $L^*_{0, \lambda}$ associated to it.
Classically
\begin{equation}L_{0,\lambda}u(x) = \int_{\mathbb{R}^d}K_\lambda(x,y)u(y)dy \end{equation}
\begin{equation}L^*_{0, \lambda}u(x) = \int_{\mathbb{R}^d}\overline{K}_\lambda(y,x)u(y)dy\end{equation}
\begin{lemma}If $u$ is smooth whose all derivative belong to $L^2$, then $(L_0-L_{0,\lambda})u$ tends to zero as well as all his derivatives and in $L^2$ 
when $\lambda \rightarrow \infty$.
 The same holds for  $(L^*_0-L^*_{0,\lambda})u$.
\end{lemma}
{\bf{Proof:}}$(L_0-L_{0,\lambda})u$ is given by the oscillatory integral
\begin{equation}\int\int_{\mathbb{R}\times \mathbb{R}^d}\exp[2\pi i<x-y\vert \xi>(1-\phi(\xi/\lambda))a(x, \xi)u(y)dyd\xi\end{equation}
The result holds by integrating by parts in $y$. Analog statement work for $(L^*_0-L^*_{0,\lambda})u$. $\diamondsuit$

{\bf{Proof of the It\^o formula:}}We put
\begin{equation}L_\lambda = L^*_{0,\lambda}L_{0, \lambda}\end{equation}
$L_\lambda$ is a continuous operator acting on bounded continuous function on $\mathbb{R}^d$ endowed with its uniform norm. The same is true for its It\^o
transform $\hat{L}_\lambda$. Therefore $L_\lambda$ generates a semi-group $P_{\lambda, t}$ on bounded continuous functions on $\mathbb{R}^d$. 
$\hat{L}_\lambda$ generates a semi-group $\hat{P}_{\lambda,t}$on bounded continuous functions on $\mathbb{R}^d \times \mathbb{R}$.
 Moreover if $u$ and $v$ are bounded continuous,
\begin{equation}P_{\lambda,t}u = \sum 1/n!L_\lambda^nu \end{equation}
and
\begin{equation}\hat{P}_{\lambda,t}v = \sum 1/n!\hat{L}_\lambda^nv \end{equation}
But
\begin{equation}L^n_\lambda\hat{v}(x) = (\hat{L}^n_\lambda v)(x, f(x))\end{equation}
Therefore
\begin{equation}P_{\lambda,t}\hat{v}(x) = (\hat{P}_{\lambda,t}v)(x,f(x))\end{equation}
But $(\hat{P}_{\lambda,t}-\hat{P}_t)(v)$ is solution of the parabolic equation
\begin{equation}-d/dt v_t = \hat{L}_\lambda v_t + (\hat{L}_{\lambda,t}-\hat{L})\hat{P}_tv\end{equation}
with initial condition 0. The result arises from the two previous lemma, by the method of variation of constants since $\hat{P}_{\lambda,t}$
 is a semi-group of contraction on $L^2(dx\otimes dy)$. This shows that for $\lambda \rightarrow \infty$
\begin{equation}\hat{P}_{\lambda,t}v \rightarrow \hat{P}_tv\end{equation}
in $L^2(dx\otimes dy)$. Similarly, in $L^2(dx)$
\begin{equation}P_{\lambda,t}\hat{v} \rightarrow P_t\hat{v}\end{equation}
We remark that $\hat{L}_\lambda$ commute with $\overline{L}_1$. Therefore
\begin{equation}(\overline{L}_1^k) (\hat{P}_{\lambda,t}-\hat{P}_t)v= (\hat{P}_{\lambda,t}-\hat{P}_t)(\overline{L}_1^kv)\end{equation}
By a similar argument to the proof of lemma (4), we can show that the convergence in (24) and (25) works for the uniform topology and not in $L^2$ only.
 This shows the result.$\diamondsuit$


\begin{thebibliography}{99}

\bibitem{1} C. Dellacherie, P.A. Meyer: Probabilit\'es et potentiel (II). Th\'eorie des  martingales. Hermann. Paris (1980).

\bibitem{2}J. Dieudonn\'e: El\'ements d'analyse VII. Gauthier-Villars. Paris (1977).

\bibitem{3}P. Gilkey: Invariance theory, the heat equation and the Atiyah-Singer theorem. Sd edition. CRC Press, Boca Raton (1995).

\bibitem{4}L. Hoermander: The analysis of linear partial differential operators (III). Springer, Heidelberg (1984)

\bibitem{5}N. Jacob: Pseudo differential operators. Markov processes (II). Generators and their potential theory. Imperial College Press. London (2002).

\bibitem{6}R. L\'eandre: Malliavin Calculus of Bismut type without probability. In "Festchrift in honour of K. Sinha".A.M. Boutet de Monvel and al eds,
Proc. Indian. Acad. Sci (Math. Sci), 116, 2006, 507-518. {\it{arXiv:0707.2143v1[math.PR]}}

\bibitem{7}R. L\'eandre: Varadhan estimates without probability: lower bounds. In "Mathematical methods in engineerings" (Ankara), D. Baleanu and al eds.
Springer, Heidelberg, 2007, 205-217.

\bibitem{8}R. L\'eandre: Positivity theorem in semi-group theory.  Mathematische Zeitschrift 258, 2008, 893-914.

\bibitem{9}R. L\'eandre: Applications of the Malliavin Calculus of Bismut type without probability.In "Simulation,
Modelling and Optimization" (Lisboa), A. M. Madureira C.D. 2006, pp. 559-564.
WSEAS transactions on mathematics 5, 2006, 1205-1211.

\bibitem{10}R. L\'eandre: The division method in semi-group theory. In "Applied mathematics" (Dallas) K. Psarris edt. W.S.E.A.S. press, 
Athens, 2007, 7-11.

\bibitem{11}R. L\'eandre: Leading term of a hypoelliptic heat-kernel. WSEAS Transactions on mathematics 6, 2007, 755-763.

\bibitem{12}R. L\'eandre: Girsanov transformation for Poisson processes in semi-group theory. In "Num. Ana. Applied. Mathematics."(Corfu) T. Simos edt.
A.I.P. Proceedings 936, 2007, 336-339.

\bibitem{13}R. L\'eandre: Malliavin Calculus of Bismut type for Poisson processes without probability.
 In "Fractional order systems". J. Sabatier and al eds. Jour. Eur. systemes Automatis\'es. 42, 2008, 715-733.

\bibitem{14}R. L\'eandre: Wentzel-Freidlin estimates in semi-group theory. in "Control, Automation, Robotics and Vision" (I.E.E.E.),
(Hanoi), Yeng Chai Soh edt, C.D., I.E.E.E.  2008, 2233-2236.

\bibitem{15}R. L\'eandre: Varadhan estimates in semi-group theory: upper bound.  "Applied computing conference" (Istanbul).M. Garcia-Planas and al eds.
W.S.E.A.S. press, 2008, 77-80.

\bibitem{16}R. L\'eandre: Varadhan estimates without probability: upper bound. WSEAS transactions on mathematics 7, 2008, 244-253. 

\bibitem {17}R. L\'eandre: Malliavin Calculus of Bismut type in semi-group theory.  Far East Journal of Mathematical Sciences 30, 2008, 1-26.

\bibitem{18}R. L\'eandre: Large deviations estimates in semi-group theory. In "Num. Ana. Applied. Mathematics" (Kos), T. Simos edt.
A.I.P. Proceedings 1048, 2008, 351-355.

\bibitem{19}R. L\'eandre: It\^o-Stratonovitch formula for a four order operator on a torus. In "Non-Euclidean Geometry and its applications" (Debrecen), S. Nagy edt, Acta Physica Debrecina 42, 
2008, 133-137.

\bibitem{20}R. L\'eandre: Regularity of a degenerated convolution semi-group without to use the Poisson process.
 To appear in "Nonlinear Science and Complexity"
(Porto) 
M. Silva and al eds. Mittag Leffler Preprint. Fall 2007. S.P.D.E. 10.

\bibitem{21}R. L\'eandre: It\^o-Stratonovitch formula for the Schroedinger equation associated to a big order operator on a torus. in "Fractional order 
differentiation" (Ankara), G. Zaslavsky and al eds,   Physica Scripta T 136, 2009, 014028.

\bibitem{22}R. L\'eandre: It\^o-Stratonovitch formula for the wave equation on a torus.  In "Computations of stochastic systems". M.A. El-Tawil edt.
Trans. Comp. Sciences VII. L.N.C.S. 5890, 2010,  68--75.

\bibitem{23}R. L\'eandre: It\^o formula for an integro differential operator without an associated stochastic process. To appear in
"ISAAC 2009" (London), J. Wirth edt.

\bibitem{24}R. L\'eandre: Wentzel-Freidlin estimates for jump processes in semi-group theory: lower bound. In "Int.Conf.Dif.Geometry.Dynamical. Systems"
(Bucuresti), V. Balan and al eds. B.S.G. Proceedings 17, 2010, 107-113.

\bibitem{25}R. L\'eandre: Wentzel-Freidlin estimates for jump processes: upper bound. To appear in "Worldcomp 10"(Las-Vegas) H. Arabnia edt

\bibitem{26}K Yosida: Functional analysis. Springer, Heidelberg, 1977.









 




 



\end{thebibliography}
\end{document}